\documentclass[oneside,leqno,12pt]{article}

\usepackage{latexsym}
\usepackage{amssymb}
\usepackage{amsmath}
\usepackage{mathrsfs}
\usepackage{graphicx}
\usepackage{verbatim}
\usepackage[all]{xy}

\newtheorem{theo}{Theorem}[section]

\newtheorem{lemma}[theo]{Lemma}

\newtheorem{corollary}[theo]{Corollary}
\newtheorem{prop}[theo]{Proposition}

\newenvironment{proof}{ \noindent \textbf{Proof.}}{\hfill $\Box$}

\def\ga{\mathfrak{a}}
\def\gb{\mathfrak{b}}
\def\gc{\mathfrak{c}}

\def\gf{\mathfrak{f}}
\def\gg{\mathfrak{g}}

\def\gk{\mathfrak{k}}
\def\gm{\mathfrak{m}}

\def\gs{\mathfrak{s}}
\def\gt{\mathfrak{t}}

\def\C{\mathbb{C}}

\def\Z{\mathbb{Z}}

\def\cA{\mathcal{A}}

\DeclareMathOperator{\Hom}{Hom}

\DeclareMathOperator{\Span}{Span}

\DeclareMathOperator{\codim}{codim}

\newcommand {\CC} {\mathbb {C}}

\newcommand {\NN} {\mathbb{Z}_{>0}}

\newcommand {\rk} {\mathrm{rk~}}

\renewcommand{\ss}{\mathfrak{s}}

\newcommand{\kk}{\mathfrak{k}}

\newcommand{\mm}{\mathfrak{m}}

\renewcommand{\gg}{\mathfrak{g}}

\renewcommand {\phi} {\varphi}

\newcommand{\triplebrace}[6]{\left\{\begin{array}{ll}#1 & #2 \\ #3 & #4 \\ #5 & #6 \end{array}\right.}
\newcommand{\refth}[1]{Theorem \ref{#1}}

\newcommand{\refcor}[1]{Corollary \ref{#1}}
\newcommand{\refprop}[1]{Proposition \ref{#1}}

\newcommand{\refeq}[1]{(\ref{#1})}

\newcommand{\Stab}{\mathrm{Stab}}

\renewcommand{\sp}{\mathrm{sp}}

\newcommand{\gl}{\mathrm{gl}}
\renewcommand{\sl}{\mathrm{sl}}
\newcommand{\so}{\mathrm{so}}

\newcommand{\inj}{\displaystyle \lim_\to}

\def\cplus{\hbox{$\subset${\raise0.3ex\hbox{\kern -0.55em ${\scriptscriptstyle +}$}}\ }}
\newcommand{\nsimeq}{\hbox{$\simeq${\raise0.3ex\hbox{\kern -0.55em ${\scriptscriptstyle /}$}}\ }}

\def\clplus{\hbox{$\subset${\raise0.3ex\hbox{\kern -0.55em ${\scriptscriptstyle +}$}}\ }}
\def\crplus{\hbox{$\supset${\raise1.05pt\hbox{\kern -0.55em ${\scriptscriptstyle +}$}}\ }}

\title{Locally semisimple and maximal subalgebras of the finitary Lie algebras $\gl(\infty)$, $\sl(\infty)$, $\so(\infty)$, and $\sp(\infty)$}
\author{Ivan Dimitrov\thanks{\, Research partially supported by NSERC Discovery Grant and by FAPESP Grant 2007/03735--0} \, and 
Ivan Penkov\thanks{\, Research partially supported by FAPESP Grant 2007/54820--8}}
\date{\empty}

\begin{document}

\maketitle

\vskip-3.5in 
\hskip1.75in
\vskip3in

\begin{abstract}
We describe all locally semisimple subalgebras and all maximal subalgebras of the 
finitary Lie algebras $\gl(\infty), \sl(\infty), \so(\infty)$, and $\sp(\infty)$. For simple finite--dimensional Lie algebras these 
classes of subalgebras
have been described in the classical works of A. Malcev and E. Dynkin.

Key words (2000 MSC): 17B05 and 17B65.
\end{abstract}

\section*{Introduction}
The simple infinite--dimensional finitary Lie algebras have been classified by A. Baranov a decade ago,
see \cite{Ba33}, \cite{Ba}, and \cite{BS}, and since then 
the study of these Lie algebras $\sl(\infty)$, $\so(\infty)$, and $\sp(\infty)$, as well of the finitary Lie algebra
$\gl(\infty)$, has been underway. So far some 
notable results on the structure of the subalgebras of  $\gl(\infty)$, $\sl(\infty)$, $\so(\infty)$, and $\sp(\infty)$ concern irreducible,
Cartan, and Borel subalgebras, see \cite{LP}, \cite{BS}, \cite{NP},\cite{DPS}, \cite{DP2}, and \cite{Da}. 
The objective of the present paper is to describe the locally semisimple subalgebras of 
$\gl(\infty)$, $\sl(\infty)$, $\so(\infty)$, and $\sp(\infty)$ (up to isomorphism, as well as in terms of their action on
the natural and conatural modules) and the maximal subalgebras of 
$\gl(\infty)$, $\sl(\infty)$, $\so(\infty)$, and $\sp(\infty)$.
 Our results extend classical results of A. Malcev, \cite{M}, and E. Dynkin, \cite{Dy1}, \cite{Dy2}, to 
infinite--dimensional finitary Lie algebras and are related to some earlier results of
A. Baranov,  A. Baranov and H. Strade, and F. Leinen and O. Puglisi.

A subalgebra $\gs$ of $\gl(\infty)$, $\sl(\infty)$, $\so(\infty)$, or $\sp(\infty)$ is locally semisimple if it is
a union of  
 semisimple finite--dimensional subalgebras.  The class of locally semisimple subalgebras is the natural analogue of the class of 
semisimple subalgebras of simple finite--dimensional Lie algebras. In the absence of Weyl's semisimplicity results 
for locally finite infinite--dimensional Lie algebras, it is a priori not clear whether a locally semisimple subalgebra 
of $\gl(\infty)$, $\sl(\infty)$, $\so(\infty)$, and $\sp(\infty)$ is itself a direct sum of simple constituents, cf. Corollary  in \cite{LP}. Theorem \ref{th1} proves 
that this is true and, moreover, that each simple constituent of a locally semisimple subalgebra 
of $\gl(\infty)$, $\sl(\infty)$, $\so(\infty)$, and $\sp(\infty)$ is either finite--dimensional or is itself isomorphic to 
$\gl(\infty)$, $\sl(\infty)$, $\so(\infty)$, or $\sp(\infty)$. The latter fact has been established earlier by A. Baranov.

The method of proof of Theorem \ref{th1} allows to prove also that if $\gg=\sl(\infty)$ (respectively, $\gg=\so(\infty)$ 
or $\sp(\infty)$) and $\gg=\inj \gs_n$ is an exhaustion of $\gg$ by semisimple finite--dimensional Lie algebras, then 
there exist $n_0$ and nested simple ideals  $\gk_n$ of  $\gs_n$  for $n > n_0$, such that 
$\inj \gk_n = \gg$,  $\gk_n \cong \sl(k_n)$ (respectively, $\gk_n=\so(k_n)$ or $\sp(k_n)$),
 and the inclusion $\gk_n\subset\gk_{n+1}$ is simply induced by an inclusion of the natural $\gk_n$--modules 
 $V(\gk_n)\subset V(\gk_{n+1})$ (cf. Corollary 5.9 in \cite{Ba2}). 

We then study the natural representation $V$ of $\gg= \gl(\infty), \sl(\infty)$, $\so(\infty)$, and $\sp(\infty)$ as a module over any 
locally semisimple subalgebra $\ss$ of $\gg$ and show that 
\begin{itemize}
\item the socle filtration of $V$ has depth  at most 2;
\item the non--trivial simple direct summands of $V$ are just natural and conatural modules over infinite--dimensional
simple ideals of $\gs$, as well as 
finite--dimensional modules over finite--dimensional ideals 
of $\gs$; each non--trivial simple constituent of $V$ as module over a simple ideal of $\gs$ occurs with finite multiplicity;
\item the module $V/V'$ is trivial. 
\end{itemize}
Similar results hold for the conatural $\gg$--module $V_*$ for $\gg = \gl(\infty)$ and $\sl(\infty)$.

We conclude the paper by a description of maximal proper subalgebras of $\gg=$ $\gl(\infty)$, $\sl(\infty)$, 
$\so(\infty)$, and $\sp(\infty)$. The maximal subalgebras of $\gg=\gl(\infty)$ are $[\gg, \gg] \cong \sl(\infty)$
and the stabilizers of subspaces of $V$ or $V_*$ as follows: $W \subset V$ with $W^{\perp\perp} = W$, or
$W \subset V$, $\codim_V W = 1$ and $W^\perp = 0$, or $\tilde{W} \subset V_*$, $\codim_{V_*} \tilde{W} = 1$ and 
$\tilde{W}^\perp = 0$. The maximal subalgebras of $\sl(\infty)$ are intersections of the maximal subalgebras of
$\gg = \gl(\infty)$ with $\sl(\infty) =[\gg, \gg]$.
For $\gg=\so(\infty)$ and $\sp(\infty)$ any maximal subalgebra is the stabilizer in $\gg$ of an isotropic subspace $W \subset V$ with $W^{\perp \perp} = W$, 
or of a non--degenerate subspace $W \subset V$ with $W \oplus W^\perp = V$ (where for $\so(\infty)$, $\dim W \neq 2$ and $\dim W^\perp \neq 2$), 
or of a non--degenerate subspace $W \subset V$ of codimension 1 such that $W^\perp = 0$.

\noindent
{\bf Acknowledgment.} We are indebted to Gregg Zuckerman for his long term encouragement to study Dynkin's papers 
\cite{Dy1} and \cite{Dy2}. 
We thank A. Baranov for very enlightening discussions and Y. Bahturin for a key reference on irreducible 
subalgebras. We acknowledge the hospitality of the Mathematisches Forschungsinstitut Oberwolfach 
where this work was initiated.

\section{General preliminaries}
The ground field is $\CC$. In this paper $V$ is a fixed countable--dimensional vector space with basis $v_1,v_2,\dots$ and $V_*$ is the restricted dual of $V$, i.e. the span of the dual set $v_1^*,v_2^*,\dots$ ($v_i^*(v_j)=\delta_{ij}$). The space $V\otimes V_*$ ($\otimes$ stands throughout the paper for 
tensor product over $\C$) has an obvious structure of an associative algebra, and by definition $\gl(V, V_*)$ (or $\gl(\infty)$ for short) is the Lie algebra associated with this associative algebra. The Lie algebra $\sl(V, V_*)$ (or $\sl(\infty)$) is  the commutator algebra $[\gl(V,V_*), \gl(V,V_*)]$.
 Given a symmetric non--degenerate form $V\times V\to\CC$, we denote by $\so(V)$ (or $\so(\infty)$) the Lie subalgebra $\Lambda^2(V)\subset\sl(V,V_*)$ (the form $V\times V\to\CC$ induces an identification of $V$ with $V_*$ which allows to consider $\Lambda^2(V)$ as a subspace of $V\otimes V_*$). Similarly, given an antisymmetric non--degenerate form $V\times V\to\CC$, we denote by $\sp(V)$ (or $\sp(\infty)$) the Lie subalgebra $S^2(V)\subset\sl(V,V_*)$. In what follows $\gg$ always stands for
one of the Lie algebras $\gl(V,V_*)$, $\sl(V,V_*)$, $\so(V)$, or $\sp(V)$.

The Lie algebras $\gl(\infty)$, $\sl(\infty)$, $\so(\infty)$, and $\sp(\infty)$ are locally finite (i.e. any finite set of elements generates a finite--dimensional subalgebra) and can be defined alternatively as follows. Recall that if $\phi:\gf\to \gf'$ is an injective homomorphism of reductive finite--dimensional Lie algebras, $\phi$ is a \emph{root injection} if for some (equivalently, for any) Cartan subalgebra $\gt_\gf$ of $\gf$, there exists a Cartan subalgebra $\gt_{\gf'}$ such that $\phi(\gt_{\gf})\subset \gt_{\gf'}$ and each $\gt_{\gf}$--root space of $\gf$ is mapped under $\phi$ into a $\gt_{\gf'}$--root space of $\gf'$. It is a known result that the direct limit $\displaystyle\lim_\to \gf_n$ of any system
\[
\gf_1\to\gf_2\to\dots
\] 
of root injections of simple finite--dimensional Lie algebras is isomorphic to $\sl(\infty)$, $\so(\infty)$, or $\sp(\infty)$, see for instance \cite{DP1}.

We need to recall also two other types of injections of simple finite--dimensional Lie algebras. Let $\gf$ and $\gf'$ be classical simple Lie algebras. We call an injective homomorphism $\phi:\gf\to\gf'$ a \emph{standard injection} if the natural representation $\omega_{\gf'}$ of $\gf'$ 
decomposes as an $\gf$--module (via $\phi$) as a direct sum of one copy of a representation which is conjugated by an automorphism of $\gf$ to the natural representation $\omega_\gf$ of $\gf$, and of a trivial $\gf$--module. Any root injection of classical Lie algebras is standard, but the converse is not true: an injection $\so(2k+1)\hookrightarrow\so(2k+2)$ is standard without being a root injection. An injective homomorphism of classical Lie algebras $\phi:\gf\to\gf'$ is \emph{diagonal}  if $\omega_{\gf'}$ decomposes as an $\gf$--module as a direct sum of copies of $\omega_\gf$, of the dual module $\omega_\gf^*$, and of the 1--dimensional trivial $\gf$--module. 
This definition is a special case of a more general definition of A. Baranov, \cite{Ba2}, \cite{BZh}.

An \emph{exhaustion} $\displaystyle\lim_\to \gg_n$ of $\gg$  is a system of injections of 
finite--dimensional Lie algebras
\[
\gg_1\stackrel{\phi_1}{\to}\gg_2\stackrel{\phi_2}{\to}\dots
\] 
such that the direct limit Lie algebra $\displaystyle\lim_\to\gg_n$ is isomorphic to $\gg$. A \emph{standard exhaustion} is an exhaustion $\displaystyle\gg=\lim_\to\gg_n$ such that $\gg_n\to\gg_{n+1}$ is a standard injection of classical simple Lie algebras for all $n$. In a standard exhaustion, for large enough $n$, $\gg_n$ is of type $A$ for $\gg=\sl(\infty)$, $\gg_n$ is of type B or D for $\gg=\so(\infty)$, and $\gg_n$ is of type $C$ for $\gg\cong\sp(\infty)$. 

A subalgebra $\ss$ of $\gg$ is \emph{locally semisimple} if it admits an exhaustion $\displaystyle\ss=\lim_\to\ss_n$ by injective homomorphisms $\ss_n\to\ss_{n+1}$ of semisimple finite--dimensional Lie algebras $\ss_n$.

For $\gg\cong \gl(\infty)$ or $\sl(\infty)$  the vector spaces $V$ and $V_*$ are by definition the \emph{natural} and \emph{conatural} $\sl(\infty)$--modules. They are characterized by the following property: $V$ (respectively, $V_*$) is the only simple $\gg$--module which, for any standard exhaustion $\displaystyle\gg=\lim_\to \gg_n$, restricts to one copy of the natural (respectively, its dual) representation of $\gg_n$ plus a trivial module. For $\gg \cong \so(\infty)$ or $\sp(\infty)$, $V$ is characterized by the same property (here $V\cong V_*$  as $\gg$--modules). 

\section{Index of a subalgebra}

For a simple finite--dimensional Lie algebra $\gf$ we denote by $\langle \cdot, \cdot \rangle_\gf$ the invariant non--degenerate symmetric bilinear 
form on $\gf$ 
for which $\langle\alpha^\vee,\alpha^\vee\rangle_\gf =2$ for any long root $\alpha$ of $\gf$. 
(By convention the roots of a simply--laced Lie algebra are long.)
 If $\phi:\gf\to\gf'$ is a homomorphism of a simple Lie algebra 
$\gf$ into the simple Lie algebra $\gf'$, then $\langle x,y\rangle_\phi:=\langle \phi(x),\phi(y)\rangle_\gf$ is an invariant symmetric bilinear 
form on $\gf$. Consequently, 
\[\langle x,y\rangle_\phi=I_\gf^{\gf'}\langle x, y\rangle_\gf\]
for some scalar $I_\gf^{\gf'}$. E. Dynkin, \cite{Dy2}, calls $I_\gf^{\gf'}$ the \emph{index of $\phi$}. 
The homomorphism $\phi$ is determined (up to an automorphism of $\gf'$) by the pull--back of any nontrivial representation of $\gf'$ of minimal dimension.
Such a representation is unique unless $\gf'$ is isomorphic to $\sl(n)$, to $D_4$,  or to $E_6$.
In the rest of the paper we fix  a non--trivial representation $\omega_{\gf'}$ of $\gf'$ of minimal dimension. 
If $\gf$ is classical, $\omega_\gf$ stands as above for the natural module.
 If $U$ is any  finite dimensional $\gf$--module, then the {\it index $I_\gf(U)$ of $U$}  is defined as  $I_\gf^{\sl(U)}$ 
where $\gf$ is mapped into $\sl(U)$ through the module $U$, see \cite{Dy2}. 
The following properties are established in \cite[$\S$ 2]{Dy2}.
\begin{prop}\label{propDyIndex} $\phantom{x}$
\begin{itemize}
\item[\rm{(i)}] $I_\gf^{\gf'}\in\Z_{\geq 0}$.
\item[\rm{(ii)}] $I_\gf^{\gf'} I_{\gf'}^{\gf''}=I_\gf^{\gf''}$.
\item[\rm{(iii)}] $I_\gf (U_1\oplus \dots \oplus U_l)=I_\gf (U_1)+\dots +I_\gf (U_l)$.
\item[\rm{(iv)}] $I_\gf (U_1\otimes \dots \otimes U_l)=\dim U_1\dots\dim U_l (\frac{1}{\dim U_1}I_\gf (U_1)+\dots +\frac{1}{\dim U_l}I_\gf (U_l))$.
\item[\rm{(v)}] If $I_\gf^{\gf'}=1$, then the root spaces of $\gf$ corresponding to long roots are mapped into root spaces of $\gf'$
corresponding to long roots. 
\end{itemize}
\end{prop}  

In particular, (ii) implies that  $I_\gf(\omega_{\gf'}) =  I_{\gf}^{\gf'} I_{\gf'}(\omega_{\gf'})$. Furthermore, a  combination of (ii)
and the information from Table 5 in \cite{Dy2} shows that 
$I_\gf^{\sp(U)} = I_\gf(U)$ and $I_\gf^{\so(U)} = \frac{1}{2} I_\gf(U)$ when $U$ admits a corresponding invariant form, see \cite{Dy2}.

We need an extension of Proposition \ref{propDyIndex}. Let $\phi:\gf\to\kk_1\oplus\dots\oplus\kk_l$ and $\eta:\kk_1\oplus\dots\oplus\kk_l\to\gf'$ be  homomorphisms of Lie algebras, where $\kk_1,\dots,\kk_l$ are simple Lie algebras. 

\begin{prop}\label{propIndex1} We have
\begin{equation}\label{eqIndex1}
I_\gf^{\gf'}=\sum_{j=1}^{l} I_\gf^{\kk_j} \, I_{\kk_j}^{\gf'},
\end{equation}
where $\gf \to \gf'$ is the homomorphism $\eta \circ \phi$, and 
the homomorphisms $\gf \to \gk_i$ and $\gk_i \to \gf'$ are determined by $\phi$ and $\eta$ in the obvious way.
\end{prop}

\begin{proof} 
Multiplying by $I_{\gf'}(\omega_{\gf'})$ we see that (\ref{eqIndex1}) is equivalent to
\begin{equation}\label{equation11}
I_\gf(\omega_{\gf'})=\sum_{j=1}^{l} I_\gf^{\kk_j} \, I_{\kk_j}(\omega_{\gf'}).
\end{equation}
In the case when $\omega_{\gf'}$ is a reducible $(\gk_1 \oplus \dots \oplus \gk_l)$--module we use Proposition \ref{propDyIndex}(iii)
to prove (\ref{equation11}) by induction on the length of $\omega_{\gf'}$.
Now assume that $\omega_{\gf'}$ is an irreducible $\kk_1\oplus\dots\oplus\kk_l$--module. 
Then $\omega_{\gf'}=U_1\otimes\dots \otimes U_l$ for some irreducible $\kk_j$--modules $U_j$. 
Note that if $U_j = \omega_{\gk_j}$  for every $j$, identity \refeq{equation11} follows from \refprop{propDyIndex}. 
Indeed, in this case $I_{\kk_j}(\omega_{\gf'})=\frac{\dim(U_1\otimes\dots\otimes U_l)}{\dim U_j} I_{\gk_j}(U_j)$ by (iii), 
and applying (iv) we obtain
$$
I_\gf(\omega_{\gf'})=\sum_j\frac{\dim(U_1\otimes\dots\otimes U_l)}{\dim U_j}I_\gf(U_j) = \sum_j \frac{I_{\gk_j}(\omega_{\gf'})}{I_{\gk_j}(U_j)} I_\gf(U_j)
= \sum_j I_\gf^{\gk_j} I_{\kk_j}(\omega_{\gf'}).
$$

To prove (\ref{equation11}) for general irreducible $\gk_j$--modules $U_j$
we consider the diagram
$$
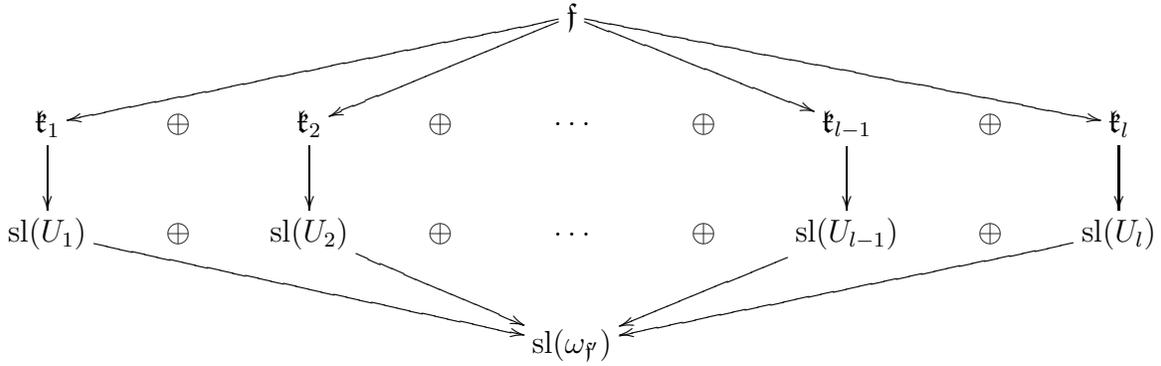
\xymatrix{
& & & & \gf\ar[lllld]\ar[lld]\ar[rrd]\ar[rrrrd] & & & & \\
\gk_1\ar[d]& \oplus & \gk_2\ar[d] & \oplus & \dots & \oplus & \gk_{l-1}\ar[d] & \oplus & \gk_l\ar[d]\\
\sl(U_1)\ar[drrrr]& \oplus & \sl(U_2)\ar[drr] & \oplus & \dots & \oplus & \sl(U_{l-1})\ar[dll] & \oplus & \sl(U_l)\ar[dllll]\\
& & & & \sl(\omega_{\gf'}) & & & &
} 
$$
This diagram enables us to first apply  (\ref{equation11}) to $\gf \to \sl(U_1) \oplus \dots \oplus \sl(U_l) \to \sl(\omega_{\gf'})$ and then
use $I_\gf^{\sl(U_j)} =  I_\gf^{\kk_j}I_{\kk_j}^{\sl(U_j)}$ to get
\[
I_\gf(\omega_{\gf'})=\sum_j I_\gf^{\sl(U_j)}I_{\sl(U_j)}(\omega_{\gf'})=\sum_j I_\gf^{\kk_j}I_{\kk_j}^{\sl(U_j)}I_{\sl(U_j)}(\omega_{\gf'})=\sum_jI_\gf^{\kk_j}I_{\kk_j}(\omega_{\gf'}).
\]
This completes the proof.
\end{proof}
\begin{prop}\label{prop3.3}
Let $\phi:\gf\to\gf'$ denote an injective homomorphism of classical simple Lie algebras. 
\begin{itemize}
\item[\rm{(i)}] Assume that $\rk \gf >4$. If $\gf'$ is not of type $B$ or $D$ and $I_\gf^{\gf'}=1$, then $\phi$ is a standard injection. Similarly,
 if $\gf$ is of type $B$ or $D$ and $I_\gf^{\gf'}=1$, then $\phi$ is a standard injection. 
\item[\rm{(ii)}] For any $n$ there exists a constant $c_n$ depending on $n$ only, such that $\rk \gf =n$ and $I_\gf^{\gf'}\leq c_n$ imply that $\phi$ is diagonal. 
Furthermore, $\lim_{n\to\infty}c_n=\infty$. 
\end{itemize}
\end{prop}
\begin{proof} 
(i) Assume first that $\gf'$ is not of type $B$ or $D$. Then $I_\gf^{\gf'} = I_\gf(\omega_{\gf'}) =1$. Proposition  \ref{propDyIndex}(iii) implies that
$\omega_{\gf'}$ considered as an $\gf$--module has exactly one non--trivial irreducible constituent $U$ with $I_\gf(U) =1$. 
We show now that $U$ is isomorphic to $\omega_\gf$ or to $\omega_\gf^*$.
Theorem 2.5 of \cite{Dy2} states that
\begin{equation} \label{equation3}
I_\gf(U) = \frac{\dim U}{\dim \gf} \langle \lambda, \lambda + 2 \rho \rangle,
\end{equation}
where 
$\langle \cdot, \cdot \rangle$ is the form induced on $\gf^*$ by $\langle \cdot, \cdot \rangle_\gf$, 
$\lambda$ is the highest weight of $U$, and $\rho$ is the half--sum of the positive roots of $\gf$. Since both $\dim U$ and 
$\langle \lambda, \lambda + 2 \rho \rangle$ are increasing functions of $\lambda$ (with respect to 
the order: $\lambda' > \lambda ''$ if $\lambda' - \lambda''$
is a non--negative combination of fundamental weights), so is $I_\gf(U)$. Table 5 in \cite{Dy2} shows that, for $\rk \gf > 4$,  a fundamental representation $U$ of $\gf$
 with $I_\gf(U) = 1$ is isomorphic to  $w_\gf$ or $\omega_\gf^*$. The monotonicity of $I_\gf(U)$ now shows that
$I_\gf(U) = 1$ implies $U \cong w_\gf$ or $U \cong \omega_\gf^*$.  Since for $\rk \gf >4$ every $\gf$--module conjugate
to $\omega_\gf$ is isomorphic to $\omega_\gf$ or $\omega_\gf^*$, $\phi$ is a standard injection.

If $\gf$ is of type $B$ or $D$, an argument similar to the one above shows that $I_\gf(\omega_{\gf'}) \geq 2$.
Consequently, formula $I_\gf^{\gf'} =  I_\gf(\omega_{\gf'}) / I_{\gf'}(\omega_{\gf'}) =1$ forces $I_\gf(\omega_{\gf'}) =2$.
Going back to the argument above we see that $I_\gf(\omega_{\gf'}) = 2$ implies that
the homomorphism $\phi$ is a standard embedding.

\noindent
(ii) Every simple Lie algebra of rank $n\geq 9$ contains a root subalgebra isomorphic to $\sl(n)$. Moreover, 
$I_{\sl(n)}^{\gf'} = \frac{I_{\sl(n)}(\omega_{\gf'})}{I_{\gf'}(\omega_{\gf'})} \geq \frac{1}{2} I_{\sl(n)} (\omega_{\gf'})$.
Hence, it is enough to show that there exist constants $d_n$ with $\lim d_n = \infty$ such that $I_{\sl(n)}(U) \geq d_n$ for any
$\sl(n)$--module $U$ which has a simple constituent not isomorphic to $\omega_{\sl(n)}$ or 
$\omega_{\sl(n)}^*$. 
To prove the existence of the constants $d_n$ we first observe that Weyl's dimension formula implies the existence a constant 
$a_1>0$, such that $\dim U \geq a_1 n^2$. Next, a direct computation gives a constant $a_2>0$, such that
$\langle\lambda,\lambda+2\rho\rangle\geq a_2 n $. Substituting these estimates into (\ref{equation3}) implies 
the existence of the constants $d_n$ with the desired properties.
\end{proof}

\begin{corollary} \label{corollary15}
Let 
$$
\gf_1 \to \gf_2 \to \dots
$$
be a system of injective homomorphisms of simple finite--dimensional Lie algebras such that $I_{\gf_n}^{\gf_{n+1}} = 1$ 
for all $n$ and $\lim (\rk \gf_n) = \infty$. Then there exists $n_0$ such that, for $n>n_0$, all homomorphisms $\gf_n \to \gf_{n+1}$ are
standard injections and all $\gf_n$ are of type $A$, or all $\gf_n$ are of type
$C$, or each $\gf_n$ is of type $B$ or $D$.
\end{corollary}

\begin{proof}
The statement follows directly from Proposition \ref{prop3.3}(ii).
\end{proof}

\section{Locally semisimple subalgebras}
\begin{theo}\label{th1}
A subalgebra $\ss\subset\gg$ is locally semisimple if and only if it is isomorphic to $\oplus_{\alpha\in A}\ss^\alpha$, where each $\ss^\alpha$ is a 
finite--dimensional simple Lie algebra or is isomorphic to $\sl(\infty)$, $\so(\infty)$, or $\sp(\infty)$, and $A$ is a finite or countable set.
\end{theo}
\begin{proof}
In one direction the statement is obvious: if $\ss\cong\oplus_{\alpha\in A}\ss^\alpha$, then by identifying $A$ with a subset of $\Z_{>0}$ and exhausting each $\ss^\alpha$ as $\displaystyle\lim_\to \ss_n^\alpha$ for some classical simple  Lie algebras $\ss_n^\alpha$, one exhausts $\ss$ via the semisimple Lie algebras $\oplus_{\alpha=1}^n\ss_n^\alpha$. 

Let now $\ss$ be locally semisimple, $\displaystyle\ss=\lim_\to\ss_n$. Write $\ss_n=\oplus_{j=1}^{l_n}\ss_n^j$, where each $\ss_n^j$ is a simple 
finite--dimensional Lie algebra. Fix a standard exhaustion $\displaystyle\gg=\lim_\to\gg_n$ of $\gg$ such that the diagram 
\begin{equation}\label{eqExhausting1}
\xymatrix{
\dots\ar[r]&\ss_n\ar[d]_{\theta_n}\ar[r]^{\phi_n}&\ss_{n+1}\ar[d]^{\theta_{n+1}}\ar[r]&\dots\ar[r]&\ss\\
\dots\ar[r]&\gg_n\ar[r]_{\psi_n}&\gg_{n+1}\ar[r]&\dots\ar[r]&\gg\\
}
\end{equation}
is commutative. In particular, $I_{\gg_n}^{\gg_{n+1}}=1$ for every $n$.

For each $1\leq j \leq l_n$ let 
\[
i_n^j:\ss_n^j\to\ss_n\quad\mathrm{and}\quad \pi_n^j:\ss_n\to\ss_n^j
\]
be the natural injection and projection respectively. Set $\theta_n^j=\theta_n\circ i_n^j:\ss_n^j \to\gg_n$ and let $\phi_n^{j,k}=\pi_{n+1}^k\circ\phi_n\circ i_n^j:\ss_n^j\to\ss_{n+1}^k$. Then $\phi_n^{i,k}$ is  a homomorphism of simple Lie algebras. Set also  
$$
\alpha_n^j:=I_{\ss_n^j}^{\gg_n}, \quad \quad \quad \beta_n^{j,k}:=I_{\ss_n^j}^{\ss_{n+1}^k}.
$$
By \refprop{propIndex1} we have
\begin{equation}\label{eqth1}
\alpha_n^j=\sum_{k=1}^{l_{n+1}}\beta_n^{j,k}\alpha_{n+1}^k.
\end{equation}
We now assign an oriented graph $\Gamma$ (a Bratteli diagram) to the direct system $\{\ss_n\}$. The vertices of $\Gamma$ are the pairs $(n,j)$ with $1\leq j\leq l_n$. A vertex $(n,j)$ has \emph{level} $n$. An arrow points from $(n,j)$ to $(n+1,k)$ if and only if $\phi_n^{j,k}$ is not trivial. 
A \emph{path} $\gamma$ in $\Gamma$ is a sequence of vertices $(n,j_n),(n+1,j_{n+1}),\dots,(m,j_m)$ such that, for every $i$ with $n\leq i\leq m-1$,  an arrow points from $(i,j_i)$ to  $(i+1,j_{i+1})$. We label the vertices and arrows of $\Gamma$ as follows: the vertex $(n,j)$ is labeled by $\alpha_n^j$ and the arrow from $(n,j)$ to $(n+1,k)$ is labeled by $\beta_n^{j,k}$. For the path $\gamma$ above we set $\gamma(i):=j_i$ for $n\leq i\leq m$ and define  $\beta(\gamma)$ as the product $\beta_n^{j_n,j_{n+1}}\beta_{n+1}^{j_{n+1},j_{n+2}} \dots\beta_{m-1}^{j_{m-1},j_m}$  of the labels of all arrows of $\gamma$. Formula \refeq{eqth1} generalizes to 
\begin{equation}
\label{eqth1-2}\alpha_n^j=\sum_\gamma\beta(\gamma)\alpha_m^{\gamma(m)},
\end{equation}
where the summation is over all paths starting at $(n,j)$ and ending at $(m,k)$ for some $1\leq k\leq l_k$. 

For each vertex $(n,j)$, let $\Gamma(n,j)$ denote the full subgraph of $\Gamma$ whose vertices appear in paths starting at $(n,j)$.
Let  $a_m(n,j)$ be the sum of the labels of all vertices of $\Gamma(n,j)$ of level $m$, i.e.
$a_m(n,j):=\sum_{(m,k)\in \Gamma(n,j)}\alpha_m^k$. Then
\begin{equation}
\label{eqth1-3}a_m(n,j)=\sum_{(m,k)\in\Gamma(n,j)}\alpha_m^k=\sum_{(m+1,t)\in\Gamma(n,j)}\left(\sum_{(m,k)\in\Gamma(n,j)} \beta_m^{k,t}\right)\alpha_{m+1}^t\geq a_{m+1}(n,j).
\end{equation}
This implies that the sequence $\{a_m(n,j)\}$ stabilizes, i.e. $a_m(n,j)=a(n,j)$ for $m$ large enough. Furthermore, \refeq{eqth1-3} shows that if $a_m(n,j)=a_{m+1}(n,j)=a(n,j)$, then each vertex  of $\Gamma(n,j)$ of level $m$ points to exactly one vertex of $\Gamma(n,j)$ of level $(m+1)$.
 In other words, the graph $\Gamma(n,j)$ is nothing but several disjoint strings from some level on.
More precisely, there exist $m_0$ and $t$ such that, for $m \geq m_0$, $\Gamma(n,j)$ has exactly $t$ vertices $(m,j_{m,1}),\dots,(m,j_{m,t})$ of level $m$
and the arrows pointing from vertices of level $m$ to vertices of level $m+1$, after a possible relabeling of the vertices of level $m+1$,
 are \[\begin{array}{lll}(m,j_{m,1})&\to& (m+1,j_{m+1,1})\\&\vdots& \\ (m,j_{m,t})&\to &(m+1,j_{m+1,t}).\end{array}\] 
 Finally, formula \refeq{eqth1-3} implies $\beta_m^{j_{m,i},j_{m+1,i}}=1$ for every $1\leq i\leq t$.

Let $\ss_m(n,j):=\oplus_{(m,k)\in\Gamma(n,j)}\ss_m^k$. Clearly $\phi_m(\ss_m(n,j))\subset\ss_{m+1}(n,j)$, hence $\displaystyle \ss(n,j)= \lim_\to \ss_m(n,j)$ is a well--defined Lie subalgebra of $\gg$. 
The fact that $\Gamma(n,j)$ splits into $t$ disjoint strings for $m\geq m_0$ implies that 
\[
\ss(n,j)=\oplus_{i=1}^t\ss^i(n,j),
\]
where $\ss^i(n,j):=\displaystyle\lim_{\substack{\to\\m\geq m_0}}\ss_m^{j_{m,i}}$. The equality $\beta_m^{j_{m,i},j_{m+1,i}}=1$ implies via Corollary 
\ref{corollary15}  that $\ss^i(n,j)$ is a 
finite--dimensional simple Lie algebra or is a Lie algebra isomorphic to $\sl(\infty)$, $\so(\infty)$, or $\sp(\infty)$. 

We are now ready to construct a decomposition $\ss=\oplus_{\alpha\in A}\ss^\alpha$ as required.
Notice first that $\Gamma(n,j) \cap \Gamma(n',j')$ is either empty or consists of several disjoint strings from some level on.
Hence $\ss(n,j)$ and $\ss(n',j')$ intersect in subsums of the direct sums $\ss(n,j)=\oplus_{i=1}^t\ss^i(n,j)$ and $\ss(n',j')=\oplus_{i'=1}^{t'}\ss^{i'}(n',j')$. Consequently,
\begin{equation} \label{equation111}
\ss=\sum_{(n,j)\in\Gamma} \gs(n,j).
\end{equation}
Let $A(n,j)$ 
denote set of paths of $\Gamma(n,j)$ and let
$\sim$ be the following equivalence relation on the set $\cup_{(n,j)\in\Gamma}A(n,j)$: $a\in A(n,j)\sim a'\in A(n',j')$ if $a$ and $a'$ coincide for large enough $m$. Define
 $A:=\left(\cup_{(n,j)\in\Gamma}A(n,j)\right)/\sim$ and, for every $\alpha \in A$, set $\ss^\alpha:=\ss^i(n,j)$, where $(m,j_{m,i}),(m+1,j_{m+1,i}),\dots$ 
 is a representative of $\alpha$. Equation (\ref{equation111}) implies that $\ss=\oplus_{\alpha\in A}\ss^\alpha$ and this 
 completes the proof.
 \end{proof}

We will illustrate the results of this paper in a series of examples built on the same set--up, cf. Theorem 5.8 in \cite{Ba1}.

\noindent\textbf{Example 1.}
Set $\tilde{V} := V \oplus \C \tilde{v}$ with $\langle \tilde{v}, v_j^*\rangle = 1$ for every $j$. Both couples $V, V_*$ and $\tilde{V}, V_*$ are non--degenerately paired
and both Lie algebras $\gg = [V \otimes V_*, V \otimes V_*]$ and  $\tilde{\gg} = [\tilde{V} \otimes V_*, \tilde{V} \otimes V_*]$ are isomorphic to $\sl(\infty)$.
Any partition $\Z_{>0} = \sqcup_{\alpha \in A} I^\alpha$ defines a locally semisimple subalgebra $\gs$ of both $\gg$ and $\tilde{\gg}$ in the following way.
Set $V^\alpha := \Span \{v_j\}_{j \in I^\alpha}$, $(V^\alpha)_* := \Span \{v_j^*\}_{j \in I^\alpha}$, and 
$\gs^\alpha := \gg \cap (V^\alpha \otimes (V^\alpha)_*)$. Define $\gs$ as $\oplus_{\alpha \in A} \gs^\alpha$.
In particular, $\gg$ itself is a locally semisimple subalgebra of $\tilde{\gg}$. 

A  corollary of \refth{th1} concerns the structure of an arbitrary exhaustion of $\gg$ by semisimple Lie algebras, cf. Corollary 5.9 in \cite{Ba2}.
\begin{corollary}\label{corth1}
Let $\displaystyle\gg=\lim_\to\ss_n$, where each $\ss_n$ is semisimple. There exist $n_0$ and simple ideals 
$\kk_n\subset\ss_n$ for $n\geq n_0$, such that $\kk_n\subset\kk_{n+1}$ and $\displaystyle\gg=\lim_\to\kk_n$. 
Furthermore, the system $\{\kk_n\}$ admits a refinement $\{\gg_s\}$ with 
\[
\gg_s\cong\triplebrace{\sl(s)}{\mathrm{if~}\gg=\sl(\infty)}{\so(s)}{\mathrm{if~}\gg=\so(\infty)}{\sp(2s)}{\mathrm{if~}\gg=\sp(\infty).}
\] 
\end{corollary}
\begin{proof}
By \refth{th1} $\gg=\oplus_{\alpha\in A}\ss^\alpha$. Since $\gg$ is simple, $A$ consists of a single element, i.e. there exists $m$ such that,
for $n\geq m$, $\Gamma(n,j)$ is a single string
\[
(m,j_m),(m+1,j_{m+1}),\dots\quad .
\]
Set $\kk_n:=\ss_n^{j_n}$ for $n \geq m$. Clearly $\displaystyle\gg=\lim_\to\kk_n$. Note that, as $I_{\gk_n}^{\gk_{n+1}} =1$ for large enough $n$, 
Corollary \ref{corollary15} implies that
there exists $n_0 \geq m$ such that
all injections $\gk_n \to \gk_{n+1}$ are standard for $n \geq n_0$.  The fact that a standard exhaustion of $\gg$ admits 
 a refinement as in the statement of the corollary is  obvious.
\end{proof}

In the special case when $\gg$ is exhausted by simple Lie algebras $\gg_n$, \refcor{corth1} implies that, for large enough $n$,
all injections $\gg_n \to \gg_{n+1}$ are standard. Furthermore, by 
 Corollary \ref{corollary15}
all $\gg_n$ are of type $A$, or all $\gg_n$ are of type $C$, or each $\gg_n$ is of type $B$ or $D$.

Here is an example showing that there exist interesting exhaustions of $\sl(\infty)$ by non--reductive Lie algebras. 

\medskip

\noindent\textbf{Example 2.} We build on Example 1.
Put $V_n := \Span \{v_1 , v_2 , \dots ,v_n \} \subset V$, $\tilde{V}_n := V_n \oplus \C \tilde{v} \subset \tilde{V}$, and
$(V_n)_* := \Span \{v_1^* , v_2^* , \dots ,v_n^* \} \subset V_*$. Set also $\gg_n = \gg \cap (V_n \otimes (V_n)_*)$ and 
$\tilde{\gg}_n := \gg \cap (\tilde{V}_n \otimes (V_n)_*)$. 
Then $\C (\tilde{v} - v_1 - \dots - v_n) \otimes (V_n)_*$ is the radical of $\tilde{\gg}_n$ and $\inj \tilde{\gg}_n$ is an exhaustion 
of $\tilde{\gg}$ with non--reductive finite dimensional Lie algebras.
Note that the Levi components $\gg_n$ of $\tilde{\gg}_n$ are nested and their direct limit $\inj \gg_n$ is nothig but the proper subalgebra
$\gg$ of $\tilde{\gg}$. On the other hand, a different choice of Levi components of $\tilde{\gg}_n$ yields an exhaustion of $\tilde{\gg}$. Indeed,
the Lie algebras $\gk_n := \tilde{\gg} \cap (\tilde{V}_{n-1} \otimes (V_n)_*)$ are also nested
and their direct limit $\inj \gk_n$ is the entire Lie algebra $\tilde{\gg}$. Moreover, since $\tilde{V}_{n-1}$ and 
$(V_n)_*$ are non--degenerately paired, we have $\gk_n \cong \sl(n)$,  which means that $\gk_n$ is a Levi component of $\tilde{\gg}_n$ for
every $n$.

We conclude this section by another corollary of \refth{th1}. 
\begin{corollary}\label{thGeneralth1} 
Let $\ga$ be a Lie algebra isomorphic to a finite or countable direct sum of finite--dimensional simple Lie algebras and of 
copies of $\sl(\infty)$, $\so(\infty)$, and $\sp(\infty)$. Then a subalgebra $\ss\subset\ga$ is locally semisimple if and only if $\ss$ itself is isomorphic to a finite or countable direct sum of finite--dimensional simple Lie algebras and of 
copies of $\sl(\infty)$, $\so(\infty)$, and $\sp(\infty)$. 
\end{corollary}

\begin{proof}
Since $\ga$ admits an obvious injective homomorphism into $\sl(\infty)$, the statement follows directly from Theorem \ref{th1}.
\end{proof}

\section{$V$ and $V_*$ as modules over a locally semisimple subalgebra $\ss \subset \gg$ }
Fix    a locally semisimple subalgebra $\gs \subset \gg$. In this section we describe the structure of $V$ and $V_*$
as  $\ss$--modules. Let $\ss=\oplus_{\alpha\in A}\ss^\alpha$ where  $\ss^\alpha$ 
are the simple constituents of $\ss$ according to \refth{th1}. Set
\begin{equation*}
\begin{array}{lcl}
A^f&:=&\{\alpha\in A\,|\, \ss^\alpha {\text { is finite--dimensional}}\},\\
A^{inf}&:=&\{\alpha\in A\,|\, \ss^\alpha {\text { is infinite--dimensional}}\},\\
\ss^f&:=&\oplus_{\alpha\in A^f}\ss^\alpha.
\end{array}
\end{equation*}
We start by describing the structure of $V$ and $V_*$ as modules over $\ss^f$. 
\begin{prop}\label{prop4}
Let $W$ be an at most countable--dimensional $\ss^f$--module with the property that, for every $x\in\ss^f$, 
the image of $x$, considered as an endomorphism of $W$, is finite--dimensional. Then
\begin{itemize}
\item[\rm{(i)}] every simple $\ss^f$--submodule of $W$ is finite--dimensional;
\item[\rm{(ii)}] $W$ has non--zero socle $W'$, hence by {\rm {(i)}} $W'$ is a direct sum of simple finite--dimensional $\ss^f$--modules.
\item[\rm{(iii)}] $W/W'$ is a trivial $\ss^f$--module.
\end{itemize}
\end{prop}
\begin{proof}
The set $A^f$ is finite or countable. If $A^f$ is finite, $\ss^f$ is a finite--dimensional semisimple Lie algebra and, by the required property on $W$, the 
$\ss^f$--module $W$ is integrable.  Hence (by a well--known extension of Weyl's semisimplicity theorem to integrable modules) $W$ is semisimple and 
all of its simple constituents are finite--dimensional. 

Assume that $A^f$ is countable and put $A^f:=\{1,2,\dots\}$. Fix an exhaustion of $\ss^f$ of the form $\ss_n^f=\ss^1\oplus\dots\oplus\ss^n$, $\ss^n$ being the simple constituents of $\ss^f$. 
If $W$ is trivial there is nothing to prove. Assume that $W$ is non--trivial. Then $W$ is a non--trivial $\gs^n$--module for some $n$.
Let $W_\kappa^n$ be a non--trivial  isotypic component  of the $\gs^n$--module $W$, i.e. an isotypic component of $W$ 
corresponding to a non--trivial simple finite--dimensional  $\gs^n$--module.
The condition on $W$ implies that $W_\kappa^n$
is finite--dimensional as otherwise the image in $W$ of any root vector of $\gs^n$ would be infinite--dimensional.
Notice that $W_\kappa^n$ is actually an $\gs^f$--submodule of $W$ since $W_\kappa^n$ is $\gs^m$--stable for all $m$. 
Furthermore, as every non--trivial simple $\gs^f$--submodule $\tilde{W}$ of $W$ contains a non--trivial $\gs^n$--submodule for
some $n$, $\tilde{W}$ is necessarily contained in $W_\kappa^n$ for some $\kappa$. This proves (i) and (ii). 

To prove (iii) we observe that the socle $W'$ of $W$ is the direct sum of a trivial module and the sum of $W_\kappa^n$ as above for all $n$ and all $\kappa$.
\end{proof}

\medskip

\noindent\textbf{Example 3.} This example shows that $W$ is not necessarily semisimple as an $\ss^f$--module, 
i.e. that $W'$ does not necessarily equal $W$. In the set--up of Example 1 consider a partition of $\Z_{>0}$ into two--element subsets.
The corresponding locally semisimple subalgebra $\gs$ of $\tilde{\gg}$ is a direct sum of infinitely many copies of $\sl(2)$ and hence $\gs^f = \gs$. 
One checks immediately that  for $W = \tilde{V}$, we have $W' = V$.

As a next step we describe the $\ss^\alpha$--module structures of $V$ and $V_*$ for $\alpha\in A^{inf}$.

\begin{prop}\label{prop2.5}~
\begin{itemize}
\item[\rm{(i)}] For any $\alpha\in A^{inf}$, the socle $V_\alpha'$ of $V$ as an $\ss^\alpha$--module is isomorphic to $k_\alpha V^\alpha\oplus l_\alpha V_*^\alpha\oplus N^\alpha$, where $k_\alpha,l_\alpha\in \NN$, $V^\alpha$ and $V_*^\alpha$ are respectively the natural and conatural representation of $\ss^\alpha$ (here $l_\alpha=0$ for $\ss^\alpha\not \cong \sl(\infty)$ ) and $N^\alpha$ is a trivial $\ss^\alpha$--module of finite or countable dimension. Similarly, for $\gg\cong \gl(\infty)$ or $\sl(\infty)$, the socle $(V_*)'_\alpha$ of $V_*$ as an $\ss^\alpha$--module is isomorphic to $k_\alpha V_*^\alpha\oplus l_\alpha V^\alpha\oplus N^\alpha_*$, where $N_*^\alpha$ is a trivial $\ss^\alpha$--module of finite or countable dimension, not necessarily equal to the dimension of $N^\alpha$. 
\item[\rm{(ii)}] $V/V'_\alpha$ and $V_*/(V_*)'_\alpha$ are trivial $\ss^\alpha$--modules.
\end{itemize}
\end{prop}

\begin{proof}
Fix standard exhaustions of $\ss^\alpha$ and $\gg$ such that the diagram  
\begin{equation*}
\xymatrix{
\dots\ar[r]&\ss_{n-1}\ar[d]\ar[r]&\ss_n\ar[d]\ar[r]&\ss_{n+1}\ar[d]\ar[r]&\dots\ar[r]&\ss\\
\dots\ar[r]&\gg_{n-1}\ar[r]&\gg_n\ar[r]&\gg_{n+1}\ar[r]&\dots\ar[r]&\gg\\
}
\end{equation*}
commutes. As in the proof of \refth{th1} we see that, for large enough $n$, $I_{\ss_n^\alpha}^{\gg_n} $ is a constant, i.e. does not depend on $n$. Therefore, by \refprop{prop4} each injective homomorphism $\ss_n^\alpha\hookrightarrow\gg_n$ is diagonal injection for large $n$, i.e.  
\begin{equation}\label{eqprop2.5-1}
V(\gg_n)=k_\alpha V(\ss_n^\alpha)\oplus l_\alpha V(\ss_n^\alpha)^*\oplus N_\alpha^n,
\end{equation}
where $V(\gg_n)$ and $V(\ss_n^\alpha)$ are the natural representation of $\gg_n$ and $\ss_n^\alpha$ respectively, the superscript $^*$ stands for dual space, $k_\alpha+l_\alpha=I_{\ss_n^\alpha}^{\gg_n}$, and $N_\alpha^n$ is a trivial $\ss_n^\alpha$--module. Furthermore
\begin{equation}\label{eqprop2.5-2}
V(\gg_n)^*=k_\alpha V(\ss_n^\alpha)^*\oplus l_\alpha V(\ss_n^\alpha)\oplus N_\alpha^n.
\end{equation}
Since $\Hom_{\ss_n^\alpha}(V(\ss_n^\alpha), V(\ss_{n+1}^\alpha)^*)=$ ~
$\Hom_{\ss_n^\alpha}(V(\ss_n^\alpha),N_\alpha^{n+1} )=$  ~
$\Hom_{\ss_n^\alpha}(V(\ss_n^\alpha)^*, V(\ss_{n+1}^\alpha))$ ~
$=\Hom_{\ss_n^\alpha}(V(\ss_n^\alpha)^*,N_\alpha^{n+1} )=0$ and $\dim \Hom_{\ss_n^\alpha}(V(\ss_n^\alpha), V(\ss_{n+1}^\alpha))=\dim\Hom_{\ss_n^\alpha}(V(\ss_n^\alpha)^*, V(\ss_{n+1}^\alpha)^*)=1 $, the fact that $V = \inj V(\gg_n)$ and $V_* = \inj V(\gg_n)^*$ implies $\dim \Hom_{\ss^\alpha}(V^\alpha, V)=k_\alpha,$ $\dim \Hom_{\ss^\alpha}(V_*^\alpha, V_*)=l_\alpha$. Therefore $k_\alpha V^\alpha\oplus l_\alpha V_*^\alpha\subset V_\alpha'$, $k_\alpha V^\alpha_*\oplus l_\alpha V^\alpha\subset (V_*)_\alpha'$. Moreover, it follows immediately from \refeq{eqprop2.5-1} and \refeq{eqprop2.5-2} that both $V_\alpha'$ and $(V_*)_\alpha'$ can only have simple constituents isomorphic to $V^\alpha$, $V^\alpha_*$ and to the 1-dimensional trivial module. This completes the proof of (i). 

Claim (ii) follows directly from (i) and from \refeq{eqprop2.5-1} and \refeq{eqprop2.5-2}.
\end{proof}

\textbf{Example 4.} This example shows that the socle of the natural representation considered as an $\gs^\alpha$--module can
also be a proper subspace. In the notations of Example 1 we can choose the subalgebra $\gs^\alpha$ of $\tilde{\gg}$ to be $\gg$.
Then $\tilde{V}' = V$ is a proper subspace of $\tilde{V}$.
Note also that the dimensions of the trivial $\gs^\alpha$--modules $N^\alpha$ and $N_*^\alpha$ are different in this case. Indeed,
$\dim N^\alpha = 1$ while $\dim N_*^\alpha = 0$.

Put now $\tilde A:=A^{inf}\cup \{f\}$ and, for every $\alpha \in \tilde{A}$,  let $V(\alpha)$ and $V_*(\alpha)$ 
denote the sum of all non--trivial simple $\gs^\alpha$--submodules of $V$ and $V_*$ respectively.

\begin{prop}\label{prop2.6}
The sums  $\sum_{\alpha\in\tilde A}V(\alpha)$ and $\sum_{\alpha\in\tilde A}V_*(\alpha)$
are direct in $V$ and $V_*$ respectively. Each $\ss^\alpha$ acts trivially on $V(\beta)$ and $V_*(\beta)$ for $\beta \neq \alpha$. Furthermore,
$V/(\oplus_{\alpha\in\tilde A}V(\alpha))$ and 
$V_*/(\oplus_{\alpha\in\tilde A}V_*(\alpha))$ are  trivial $\ss$--modules.
\end{prop}
\begin{proof} We will prove the proposition for $V$ as the statements for $V_*$ are analogous.
Let $\alpha, \beta \in A^{inf}$ and let   $\gs^\alpha = \inj \gs^\alpha_n$ and $\gs^\beta = \inj \gs^\beta_n$ be standard exhaustions.
Assume that the action of $\ss^\alpha$ on $V(\beta)$ is non--trivial. 
Then, for some $i$, $V$ will have simple $\ss_i^\alpha\oplus\ss^\beta_n$--submodules of the form 
$V_i^\alpha\otimes M_n^\beta$ or $(V_*^\alpha)_i\otimes M_n^\beta$ for some $\ss_\beta^n$--modules $M_n^\beta$ of unbounded dimension when $n\to\infty$. 
This would imply that the multiplicity of $V_i^\alpha$ or $(V_*^\alpha)_i$ in $V$ is infinite, which is a contradiction. The case when $\alpha=f$ or $\beta=f$ is dealt with in a similar way. 

The fact that $V/\oplus_{\alpha\in\tilde A}V(\alpha)$ is a trivial $\ss$--module is obvious.
\end{proof}

In this way we have proved the following theorem. 

\begin{theo}
The socle $V'$ of $V$ (respectively, $(V_*)'$ of $V_*$) considered as an $\ss$--module 
 is isomorphic to the direct sum of all non--trivial $\ss^\alpha$--submodules $V(\alpha)$ (respectively, $V_*(\alpha)$) 
 of $V$ (respectively, $V_*$), described in 
 Propositions {\rm {\ref{prop2.5}}} and {\rm {\ref{prop2.6}}} plus a possible trivial $\ss$--submodule. The quotients $V/V'$ 
and $V_*/(V_*)'$ are trivial $\ss$--modules.  
\end{theo}
\begin{proof}
By \refprop{prop2.6}, for each $\alpha \in A^{inf}$, the modules $V(\alpha) \subset V$ and $V_*(\alpha) \subset V_*$ are semisimple $\ss$--submodules of finite length.
Moreover, the modules $V(f) \subset V$ and $V_*(f) \subset V_*$ are semisimple $\ss$--submodules with finite--dimensional simple constituents.
By \refprop{prop2.6}, the quotients  $V/\oplus_{\alpha \in \tilde{A}} V(\alpha)$ and  $V_*/\oplus_{\alpha \in \tilde{A}} V_*(\alpha)$
are trivial $\ss$--modules, and the statement follows. 
\end{proof}

Note that to any locally semisimple subalgebra $\gs \subset \gg$ we can assign some "standard invariants". These are 
the isomorphism classes of $V(f)$ and $V_*(f)$ as $\gs^f$--modules, the pairs of numbers 
$\{k_\alpha, l_\alpha\}_{\alpha \in A^{inf}}$, and the dimensions  $\{ \dim N^J, \dim N_*^J, \dim V/V_J', \dim V_*/(V_*)_J'\}_{J \subset A^{inf}}$, 
where $N^J := \cap_{\beta \in J} N^\beta$, $N_*^J := \cap_{\beta \in J} N_*^\beta$, and $V_J'$ and $(V_*)_J'$ are the respective socles of $V$ and $V_*$ considered as $(\oplus_{\beta \in J} \gs^\beta)$--modules.
Clearly, these invariants are preserved when conjugating by elements of the group $GL(V, V_*)$ of all automorphisms of $V$ under which $V_*$ is stable
(respectively, all automorphisms of $V$ preserving the non--degenerate form $V \times V \to \C$ for $\gg = \so(V)$ or $\sp(V)$). In a similar way, when
$\gs$ is replaced by a maximal toral subalgebra, it is shown in \cite{DPS} that the analogous invariants are only rather rough invariants of the $GL(V, V_*)$--conjugacy
classes of maximal toral subalgebras. The $GL(V, V_*)$--conjugacy classes of locally semisimple subalgebras $\gs \subset \gg$ with fixed "standard invariants" 
as above remain to be studied.

\section{Maximal subalgebras}
\begin{theo}$\phantom{x}$ \newline
Let $\mm\subset\gg$ be a proper subalgebra. 
\begin{itemize}
\item[\rm{(i)}] If $\gg=\gl(V, V_*)$, then $\gm$ is maximal if and only if  one of the following three mutually exclusive statements holds:
\begin{itemize}
\item[\rm{(ia)}] $\mm=[\gg,\gg]= \sl(V, V_*)$;
\item[\rm{(ib)}] $\mm=\Stab_\gg W$ or $\mm=\Stab_\gg \tilde W$, where $W\subset V$ (respectively, $\tilde W\subset V_*$) is a subspace with the properties $\codim_VW=1$, $W^\perp=0$ (respectively, $\codim_{V_*}\tilde W=1$, $\tilde W^\perp=0$);
in this case $\gm \cong \gl(\infty)$;
\item[\rm{(ic)}] $\mm=\Stab_\gg W=\Stab_\gg W^\perp$, where $W\subset V$ is a proper subspace with $W^{\perp\perp}=W$.
\end{itemize}
\item[\rm{(ii)}] If $\gg=\sl(V,V_*)$, then $\gm$ is maximal if and only if  one of the following three mutually exclusive statements holds:
\begin{itemize}
\item[\rm{(iia)}] $\mm = \so(V)$ or $\gm = \sp(V)$ for an appropriate non--degenerate symmetric or skew--symmetric form on $V$;
\item[\rm{(iib)}] $\mm=\Stab_\gg W$ or $\mm=\Stab_\gg \tilde W$, where $W\subset V$ (respectively, $\tilde W\subset V_*$) is a subspace with the properties $\codim_VW=1$, $W^\perp=0$ (respectively, $\codim_{V_*}\tilde W=1$, $\tilde W^\perp=0$);
in this case $\gm \cong \sl(\infty)$;
\item[\rm{(iic)}] $\mm=\Stab_\gg W=\Stab_\gg W^\perp$, where $W\subset V$ is a proper subspace with $W^{\perp\perp}=W$.
\end{itemize}
\item[\rm{(iii)}]  If $\gg=\so(V)$ or $\gg = \sp(V)$, then $\gm$ is maximal if and only if  $\mm=\Stab_\gg W$ for some subspace $W \subset V$ satisfying one
 of the following three mutually exclusive conditions:
\begin{itemize}
\item[\rm{(iiia)}] $W$ is non--degenerate such that  $W \oplus W^\perp = V$ and $\dim W \neq 2$, $\dim W^\perp \neq 2$ for $\gg = \so(V)$; 
in this case $\gm = \so(W) \oplus \so(W^\perp)$ when $\gg = \so(V)$, and 
 $\gm = \sp(W) \oplus \sp(W^\perp)$ when $\gg = \sp(V)$;
\item[\rm{(iiib)}] $W$ is non--degenerate such that $W^\perp = 0$ and $\codim_V W = 1$; 
in this case $\gm = \so(W)$ when $\gg = \so(V)$, and 
$\gm = \sp(W)$ when $\gg = \sp(V)$;
\item[\rm{(iiic)}]  $W$ is isotropic with $W^{\perp \perp} = W$.
\end{itemize}
\end{itemize}
The space $W$ (respectively, $\tilde{W}$) is unique in cases {\rm {(ib) and (iib)}}; the space $W$ is unique in cases {\rm {(ic), (iic), (iiib), and (iiic)}};
the pair $(W, W^\perp)$ is unique in case {\rm {(iiia)}}.
\end{theo}
\begin{proof}
Let $\gg=\gl(V, V_*)$ and let $\gm$ be maximal. 
If both $V$ and $V_*$ are irreducible $\gm$--modules, then $\gm = [\gg, \gg]$. This follows from the description of irreducible subalgebras of $\gg$ 
given in Theorem 1.3 in \cite{BS}.
Let $V$ be a reducible $\gm$--module. 
Then $\mm \subset \Stab_\gg W$ for some proper subspace $W\subset V$. 
 Since $V$ is an irreducible $\gg$--module, 
$\Stab_\gg W$ is a proper subalgebra of $\gg$. Therefore the maximality of $\mm$ yields 
$\mm=\Stab_\gg W$. If $W^{\perp\perp}=W$, we are in case (ic). If the inclusion 
$W\subset W^{\perp\perp}$ is proper, then the inclusion $\Stab_\gg W\subset\Stab_\gg W^{\perp\perp}$ is 
also proper since $W^{\perp\perp} \otimes V_* \subset \Stab_\gg W^{\perp \perp}$ and 
$W^{\perp\perp} \otimes V_* \not \subset \Stab_\gg W$. Hence we have a contradiction unless 
$\Stab_\gg W^{\perp\perp}=\gg$. In the latter case $W$ must have codimension 1 in 
$V$ as otherwise $\Stab_\gg W$ again would not be maximal.
Moreover, $\Stab_\gg W = W \otimes V_*$ and, as $W$ and $V_*$ are non--degenerately paired, $\gm = \Stab_\gg W \cong \gl(\infty)$.

Finally, if $V_*$ is a reducible $\gm$--module and $V$ is an irreducible $\gm$--module then $\gm = V \otimes \tilde{W}$ for
a subspace $\tilde{W} \subset V_*$ as in (ib). 
This proves (i) in one direction. 

For the other direction, one needs to show that if $W$ (respectively, $\tilde{W}$) is a subspace as in (ib) or (ic), $\Stab_\gg W$ (respectively, 
$\Stab_\gg \tilde{W}$) is a maximal subalgebra. In case (ic) this 
follows from the observation that  $\Stab_\gg W = W \otimes V_* + V \otimes W^\perp$
which shows that $\gg/\Stab_\gg W \cong (V/W) \otimes (V_*/W^\perp)$ is an irreducible $\Stab_\gg W$--module. In case (ib)  
$\Stab_\gg W = W \otimes V_*$ (respectively, $\Stab_\gg \tilde{W} = V \otimes \tilde{W}$), hence 
$\gg/\Stab_\gg W \cong V_*$ (respectively, $\gg/\Stab_\gg \tilde{W} \cong V$) 
is an irreducible $\Stab_\gg W$--module. The proof of (i) is now complete.

Claim (ii) is proved in the same way.

Let $\gg = \so(V)$ or $\gg  =\sp(V)$ and let $\gm$ be maximal.
Then $V$ must be a reducible $\gm$--module by Theorem 1.3 in \cite{BS}. If $W$ is a proper $\gm$--submodule of
$V$, then $\gm$ stabilizes $W^{\perp \perp}$ as well. If $W^{\perp \perp} = V$, i.e. $W^\perp = 0$, the inclusion
$\Stab_\gg W \subset \Stab_\gg W^\sharp$ is proper whenever $W$ is a proper subspace of $W^\sharp$.  
The
maximality of $\gm$ implies then $\codim_V W = 1$ and we are in case (iiic). If $W^{\perp \perp}$
is a proper subspace of $V$, the inclusions $\gm \subset \Stab_\gg W \subset \Stab_\gg W^{\perp \perp}$ and the maximality
of $\gm$ imply that $\gm = \Stab_\gg W^{\perp \perp}$. Noting that $(W^{\perp \perp})^{\perp \perp} = W^{\perp \perp}$ we may replace
$W$ by $W^{\perp \perp}$ and for the rest of the proof assume that $\gm = \Stab_\gg W$, where $W^{\perp \perp} = W$. 

If $W$ is isotropic or $W^\perp$ is isotropic, then $\Stab_\gg W=\Stab_\gg W^\perp$ and we are in case (iiic). 
 
If $W\cap W^\perp$ is a proper subspace both of $W$ and $W^\perp$, $W\cap W^\perp$ 
is an isotropic space. The inclusion $\mm\subset\Stab_\gg(W\cap W^\perp)$  implies 
$\mm=\Stab_\gg(W\cap W^\perp)$, and again we are in case (iiic) as $(W \cap W^\perp)^{\perp \perp} = W \cap W^\perp$.
Assume $W\cap W^\perp = 0$. 
Then $\gm \subset \Stab_\gg (W \oplus W^\perp)$. If $W \oplus W^\perp = V$ and $\dim W \neq 2$ and $\dim W^\perp \neq 2$ for $\gg = \so(V)$, 
then $\Stab_\gg W = \so(W) \oplus \so(W^\perp)$ or
 $\Stab_\gg W = \sp(W) \oplus \sp(W^\perp)$, and we 
are in case (iiia). The case when $\gg = \so(V)$ and $\dim W =2$ or $\dim W^\perp = 2$ does not occur as then $\Stab_\gg W$ is contained properly in
the stabilizer of an isotropic  subspace of $W$ or $W^\perp$ respectively. 

If the inclusion $W \oplus W^\perp \subset V$ is proper, 
 then $\Stab_\gg (W \oplus W^\perp)$ is a proper subalgebra of $\gg$ and the 
 the inclusion $\Stab_\gg W \subset
\Stab_\gg(W \oplus W^\perp)$ is also proper.  
Indeed, for $\gg = \so(V)$ we have $\Lambda^2(W \oplus W^\perp) \subset \Stab_\gg (W \oplus W^\perp)$ and 
$\Lambda^2(W \oplus W^\perp) \not \subset \Stab_\gg W$, and for 
$\gg = \sp(V)$ we have $S^2(W \oplus W^\perp) \subset \Stab_\gg (W \oplus W^\perp)$ and 
$S^2(W \oplus W^\perp) \not \subset \Stab_\gg W$.
Hence the maximality of $\gm$ implies $V = W \oplus W^\perp$, and we have proved (iii) in one direction.

We leave it to the reader to verify that, for every $W$ as in (iiia), (iiib), and (iiic), $\Stab_\gg W$  is a maximal subalgebra of $\gg$.

To prove the uniqueness of $W$ (respectively, $\tilde{W}$) or of the pair $(W, W^\perp)$ as stated,  it is enough to notice that $W$ (respectively, $\tilde{W}$)
is the unique proper $\gm$--submodule of $V$ (respectively, $V_*$) in cases (ib) and (iib); that $W$ 
is the unique proper $\gm$--submodule of $V$ in cases (ic), (iic), (iiib), and (iiic); and that $W$ are $W^\perp$ are the only proper
$\gm$--submodules of $V$ in case (iiia).
\end{proof}

Note that the subalgebra $\gg \subset \tilde{\gg}$ from Example 2  is a maximal simple subalgebra of $\tilde{\gg}$ as in (ib).
Furthermore, in all cases but (ic), (iic), and (iiic), a maximal subalgebra $\gm$ is irreducible in the sense of \cite{LP} and \cite{BS},
and in all cases but (ib),  (iib), and (iiib) $\gg$ admits a standard exhaustion $\inj \gg_n$ such that the Lie algebras $\gm \cap \gg_n$
are maximal subalgebras of $\gg_n$ for all $n$.

\vskip 2 cm

\centerline{\begin{tabular}{ll}
I.D.: & I.P.:\\
Department of Mathematics and Statistics& Jacobs University Bremen\\
Queen's University& Campus Ring 1\\
Kingston, K7L 3N6 & 28759 Bremen\\
Canada & Germany\\
{\tt dimitrov@mast.queensu.ca} & {\tt i.penkov@jacobs-university.de}
\end{tabular}}

\end{document}